\documentclass[a4paper,twoside,english]{amsart}
\usepackage[active]{srcltx}
\usepackage{amsthm}
\usepackage{amssymb}
\usepackage{esint}

\makeatletter

\newcommand{\noun}[1]{\textsc{#1}}

\theoremstyle{plain}
\newtheorem{theorem}{\bf Theorem}[section]

\makeatother

\usepackage{babel}

\begin{document}

\title[A note on poly-Bernoulli numbers and polynomials of the second kind ]{A note on poly-Bernoulli numbers and polynomials of the second kind }

\author[T. Kim, S. H. Lee, and J. J. Seo]{Taekyun Kim, Sang Hun Lee, and Jong Jin Seo}
\begin{abstract}
In this paper, we consider the poly-Bernoulli numbers and polynomials of the second kind and presents new and explicit formulae for calculating the poly-Bernoulli numbers of the second kind and the Stirling numbers of the second kind.
\end{abstract}

\thanks{*Corresponding Author: }
\maketitle

\section{Introduction}
As is well known, the Bernoulli polynomials of the second kind are defined by the generating function to be
\begin{equation}
\frac{t}{\log(1+t)}(1+t)^x=\sum_{n=0}^{\infty}b_{n} (x)\frac{t^n}{n!}, \ \textnormal{(see [5,14,16])}.
\end{equation}
When $x=0, b_{n}=b_{n}(0)$ are called the Bernoulli numbers of the second kind. The first few Bernoulli numbers $b_{n}$ of the second kind are
$b_{0}=1$, $b_{1}={1}/{2}$, $b_{2}=-{1}/{12}$, $b_{3}={1}/{24}$, $b_{4}=-{19}/{720}$, $b_{5}={3}/{160}, \cdots$.

From (1), we have

\begin{equation}
b_{n}(x)=\sum_{l=0}^{n}\binom{n}{l}b_{l}\ (x)_{n-l},
\end{equation}
where $(x)_{n}=x(x-1)\cdots(x-n+1), (n\geqq0)$. The Stirling number of the second kind is defined by

\begin{equation}\begin{split}
x^{n}=\sum_{l=0}^{n}S_{2}(n,l)(x)_{l},\ (n\geqq0).
\end{split}\end{equation}

The ordinary Bernoulli polynomials are given by

\begin{equation}\begin{split}
\frac{t}{e^{t}-1}e^{xt}=\sum_{n=0}^{\infty}B_{n}(x)\frac{t^{n}}{n!},\ \textnormal{(see [1-18])}.
\end{split}\end{equation}
When $x=0, B_{n}=B_{n}(0)$ are called the Bernoulli numbers.

It is known that the classical polylogarithmic function $Li_{k}(x)$ is given by

\begin{equation}\begin{split}
Li_{k}(x)=\sum_{n=1}^{\infty}\frac{x^{n}}{n^{k}}, (k\in \mathbb Z),\ \textnormal{(see [6,7,8])}.
\end{split}\end{equation}
For $k=1, Li_{1}(x)=\sum_{n=1}^{\infty}\frac{x^{n}}{n}=-\log(1-x)$. The Stirling number of the first kind is defined by
\begin{equation}\begin{split}
(x)_{n}=\sum_{l=0}^{n}S_{1}(n,l)x^{l}, (n\geq0),\ \textnormal{(see [15])}.
\end{split}\end{equation}

In this paper, we consider the poly-Bernoulli numbers and polynomials of the second kind and presents new and explicit formulae for calculating the poly-Bernoulli number and polynomial and the Stirling number of the second kind.

\section{poly-Bernoulli numbers and polynomials of the second kind}

For $k\in \mathbb Z$, we consider the poly-Bernoulli polynomials $b_{n}^{(k)}(x)$ of the second kind as follows:

\begin{equation}
\frac{Li_{k}(1-e^{-t})}{\log(1+t)}(1+t)^{x}=\sum_{n=0}^{\infty}b_{n}^{(k)}(x)\frac{t^{n}}{n!}.
\end{equation}
When $x=0, b_{n}^{(k)}=b_{n}^{(k)}(0)$ are called the poly-Bernoulli numbers of the second kind.

Indeed, for $k=1$, we have
\begin{equation}
\frac{Li_{k}(1-e^{-t})}{\log(1+t)}(1+t)^{x}=\frac{t}{\log(1+t)}(1+t)^{x}=\sum_{n=0}^{\infty}b_{n}(x)\frac{t^{n}}{n!}.
\end{equation}
By (7) and (8), we get
\begin{equation}
b_{n}^{(1)}(x)=b_{n} (x),\ (n\geq 0).
\end{equation}

It is known that

\begin{equation}
\frac{t(1+t)^{x-1}}{\log(1+t)}=\sum_{n=0}^{\infty}B_{n}^{(n)}(x)\frac{t^{n}}{n!},
\end{equation}
where $B_{n}^{(\alpha)}(x)$ are the Bernoulli polynomials of order $\alpha$ which are given by the generating function to be
\begin{equation*}
\left(\frac{t}{e^{t}-1}\right)^{\alpha}e^{xt}=\sum_{n=0}^{\infty}B_{n}^{(\alpha)}(x)\frac{t^{n}}{n!},\ \textnormal{(see [1-18])}.
\end{equation*}
By (1) and (10), we get
$$b_{n}(x)=B_{n}^{(n)}(x+1),\ (n\geq 0).$$
Now, we observe that

\begin{equation}\begin{split}
&\frac{Li_{k}(1-e^{-t})}{\log(1+t)}(1+t)^{x}\\
=&\sum_{n=0}^{\infty}b_{n}^{(k)}(x)\frac{t^{n}}{n!}\\
=&\frac{1}{\log(1+t)}\underbrace{\int_{0}^{t}\frac{1}{e^{x}-1}\int_{0}^{t}\frac{1}{e^{x}-1}
\cdots\frac{1}{e^{x}-1}}_{k-1\ times}\int_{0}^{t}\frac{x}{e^{x}-1}dx\cdots dx(1+t)^{x}.
\end{split}\end{equation}
Thus, by (11), we get

\begin{equation}\begin{split}
\sum_{n=0}^{\infty}b_{n}^{(2)}(x)\frac{t^{n}}{n!}&=\frac{(1+t)^{x}}{\log(1+t)}\int_{0}^{t}\frac{x}{e^{x}-1}dx \\
&=\frac{(1+t)^{x}}{\log(1+t)}\sum_{l=0}^{\infty}\frac{B_{l}}{l!}\int_{0}^{t}x^{l}dx \\
&=\left(\frac{t}{\log(1+t)}\right)(1+t)^{x}\sum_{l=0}^{\infty}\frac{B_{l}}{(l+1)}\frac{t^{l}}{l!}\\
&=\sum_{n=0}^{\infty}\left\{\sum_{l=0}^{n}\binom{n}{l}\frac{B_{l}b_{n-l} (x)}{l+1}\right\}\frac{t^{n}}{n!}.
\end{split}\end{equation}
Therefore, by (12), we obtain the following theorem.

\begin{theorem} For $n\geq 0$ we have
\begin{equation*}
b_{n}^{(2)}(x)=\sum_{l=0}^{n}\binom{n}{l}\frac{B_{l}b_{n-l} (x)}{l+1}.
\end{equation*}
\end{theorem}
From (11), we have

\begin{equation}\begin{split}
\sum_{n=0}^{\infty}b_{n}^{(k)}(x)\frac{t^{n}}{n!}=&\frac{Li_{k}(1-e^{-t})}{\log(1+t)}(1+t)^{x}\\
=&\frac{t}{\log(1+t)}\frac{Li_{k}(1-e^{-t})}{t}(1+t)^{x}.
\end{split}\end{equation}
We observe that
\begin{equation}\begin{split}
\frac{1}{t}Li_{k}(1-e^{-t})&=\frac{1}{t}\sum_{n=1}^{\infty}\frac{1}{n^{k}}(1-e^{-t})^{n}\\
&=\frac{1}{t}\sum_{n=1}^{\infty}\frac{(-1)^{n}}{n^{k}}n!\sum_{l=n}^{\infty}S_{2}(l,n)\frac{(-t)^l}{l!}\\
&=\frac{1}{t}\sum_{l=1}^{\infty}\sum_{n=1}^{l}\frac{(-1)^{n+l}}{n^{k}}n!S_{2}(l,n)\frac{t^{l}}{l!}\\
&=\sum_{l=0}^{\infty}\sum_{n=1}^{l+1}\frac{(-1)^{n+l+1}}{n^{k}}n!\frac{S_{2}(l+1,n)}{l+1}\frac{t^{l}}{l!}.
\end{split}\end{equation}
Thus, by (10) and (14), we get

\begin{equation}\begin{split}
\sum_{n=0}^{\infty}b_{n}^{(k)}(x)\frac{t^{n}}{n!}
&=\left(\sum_{m=0}^{\infty}b_{m}(x)\frac{t^{m}}{m!}\right)
\left\{\sum_{l=0}^{\infty}\left(\sum_{p=1}^{l+1}\frac{(-1)^{p+l+1}}{p^{k}}p!\frac{S_{2}(l+1,p)}{l+1}\right)\frac{t^{l}}{l!}\right\}\\
&=\sum_{n=0}^{\infty}\left\{\sum_{l=0}^{n}\binom{n}{l}\left(\sum_{p=1}^{l+1}\frac{(-1)^{p+l+1}p!}{p^{k}}\frac{S_{2}(l+1,p)}{l+1}\right)b_{n-l} {(x)}\right\}\frac{t^{n}}{n!}.
\end{split}\end{equation}
Therefore, by (15), we obtain the following theorem.

\begin{theorem} For $n\geq 0$, we have
\begin{equation*}
b_{n}^{(k)}(x)=\sum_{l=0}^{n}\binom{n}{l}\left(\sum_{p=1}^{l+1}\frac{(-1)^{p+l+1}}{p^{k}}p!\frac{S_{2}(l+1,p)}{l+1}\right)b_{n-l}(x).
\end{equation*}
\end{theorem}
By (7), we get
\begin{equation}\begin{split}
\sum_{n=0}^{\infty}\left(b_{n}^{(k)}(x+1)-b_{n}^{(k)}(x)\right)\frac{t^{n}}{n!}
=&\frac{Li_{k}(1-e^{-t})}{\log(1+t)}(1+t)^{x+1}-\frac{Li_{k}(1-e^{-t})}{\log(1+t)}(1+t)^{x}\\
=&\frac{tLi_{k}(1-e^{-t})}{\log(1+t)}(1+t)^{x}\\
=&\left(\frac{t}{\log(1+t)}(1+t)^{x}\right)Li_{k}(1-e^{-t})\\
=&\left(\sum_{l=0}^{\infty}\frac{b_{l}(x)}{l!}t^{l}\right)\left\{\sum_{p=1}^{\infty}\left(\sum_{m=1}^{p}\frac{(-1)^{m+p}m!}{m^{k}}S_{2}(p,m)\right)\right\}\frac{t^{p}}{p!} \\
\end{split}\end{equation}

\begin{equation}\begin{split}
&=\sum_{n=1}^{\infty}\left(\sum_{p=1}^{n}\sum_{m=1}^{p}\frac{(-1)^{m+p}}{m^{k}}m!S_{2}(p,m)\frac{b_{n-p} {(x)}n!}{(n-p)!p!}\right)\frac{t^{n}}{n!}\\
&=\sum_{n=1}^{\infty}\left\{\sum_{p=1}^{n}\sum_{m=1}^{p}\binom{n}{p}\frac{(-1)^{m+p}m!}{m^{k}}S_{2}(p,m)b_{n-p}(x)\right\}\frac{t^{n}}{n!}.
\end{split}\end{equation}
Therefore, by (16), we obtain the following theorem.

\begin{theorem} For $n\geq 1$, we have
\begin{equation}
b_{n}^{(k)}(x+1)-b_{n}^{(k)}(x)=\sum_{p=1}^{n}\sum_{m=1}^{p}\binom{n}{p}\frac{(-1)^{m+p}m!}{m^{k}}S_{2}(p,m)b_{n-p}(x).
\end{equation}
\end{theorem}
From (13), we have

\begin{equation}\begin{split}
\sum_{n=0}^{\infty}b_{n}^{(k)}(x+y)\frac{t^{n}}{n!}&=\left(\frac{Li_{k}(1-e^{-t})}{\log(1+t)}\right)^{k}(1+t)^{x+y} \\
&=\left(\frac{Li_{k}(1-e^{-t})}{\log(1+t)}\right)^{k}(1+t)^{x}(1+t)^{y} \\
&=\left(\sum_{l=0}^{\infty}b_{l}^{(k)}(x)\frac{t^{l}}{l!}\right)\left(\sum_{m=0}^{\infty}(y)_{m}\frac{t^{m}}{m!}\right) \\
&=\sum_{n=0}^{\infty}\left(\sum_{l=0}^{n}(y)_{l}b_{n-l}^{(k)}(x)\frac{n!}{(n-l)!l!}\right)\frac{t^{n}}{n!} \\
&=\sum_{n=0}^{\infty}\left(\sum_{l=0}^{n}\binom{n}{l}b_{n-l}^{(k)}(x)(y)_{l}\right)\frac{t^{n}}{n!}. \\
\end{split}\end{equation}
Therefore, by (17), we obtain the following theorem.

\begin{theorem} For $n\geq 0$, we have
\begin{equation*}
b_{n}^{(k)}(x+y)=\sum_{l=0}^{n}\binom{n}{l}b_{n-l}^{(k)}(x)(y)_{l}.
\end{equation*}
\end{theorem}

\bigskip

{\hskip -1pc \bf Acknowledgements}.The present Research has been conducted by the Research Grant of Kwangwoon University in 2014

\bigskip
\medskip
\bigskip
\medskip

\noindent \noun{T. Kim\\
 Department of Mathematics\\
Kwangwoon University\\
Seoul 139-701, Republic of Korea}\\
\textnormal{e-mail: tkkim@kw.ac.kr}

\vskip 1pc

\noindent \noun{S. H. Lee\\
Division of General Education\\
Kwangwoon University\\
Seoul 139-701, Republic of Korea}\\
\textnormal{e-mail: tkkim@kw.ac.kr}

\vskip 1pc

\noindent \noun{J. J. Seo\\
Department of Applied Mathematics\\
Pukyong National University\\
Pusan 698-737, Republic of Korea}\\
\textnormal{e-mail: seo2011@pknu.ac.kr}

\end{document}